\documentstyle[12pt,prl,aps,graphicx]{revtex}

\begin{document}
 \title{ \Large 
Random triangle in square: geometrical approach}
\normalsize \author{  Zakir F. Seidov, 
Dept of Physics, POB 653 Ben-Gurion University, \\84105 Beer-Sheva, Israel\\
E-mail: seidov@bgumail.bgu.ac.il}
\maketitle \begin{abstract} 
The classical problem of mean area of random triangle within the square is
solved by a simple and explicit method. Some other related problems are also
solved using Mathematica. \end{abstract}
\section{Intro}
We call our approach geometrical as instead of considering 6-fold integral
in abstract space we consider random triangle (RT) inside the  plane
rectangle when all possible cases are explicitly apparent.\\
Area of triangle with vertices p1=(x1,y1), p2=(x2,y2), p3=(x3,y3) is equal
to\\ $s=\frac 12(x1(y2-y3)+x2\,(-y1+y3)+x3\,(y1-y2)).\,\,(1)$\\
Let points p1, p2, p3 are randomly (with constant differential probability
function) distributed 
over the rectangle with sides A, B. What is the mean area of triangles
with vertices p1,p2,p3?\\
Answer is evident: zero, as any given triangle corresponds to 6 cases of
full permutation of three points at vertices of the triangle. Mean area
of this 6 triangles, as given by (1), is zero.\\
But if we take triangle as geometrical figure and if we consider an area
of such a figure as positive
value, then we must take absolute value of s in formula (1) and...
calculation of relevant integrals
become impossible even for Mathematica. So Michael Trott in his recent
brilliant paper in Mathematica Journal [1] found 496 different integrals
each over subregion with the same sign of $s$,
and then used Mathematica to solve such an enormously difficult task.
Needless to say that M.Trott's stunning skill of using Mathematica is
far out of scope of ordinary reader (as me, e.g.),
so I've spend some three weeks in searching a more simple solution. The
result is most easily get by the explicit geometrical approach.
\section{Geometrical approach}
Here we consider RT in square (=right rectangle) with side length $A$.
First observation is that due to points symmetry we may simplify problem by
considering particular relation between points. For example, as we will do
here, we may consider only case $x1<x2<x3$
{\it with due account of normalizing condition}.\\
First point (p1) may take any position inside the square, so the region of
integration over x1, y1 is $0<x1<A;0<y1<A$ at all cases considered
further.\\
Now we should discriminate two cases of relation between ordinates of 1st
and 2nd points: 1) $y2>y1$ and 2) $y2<y1$.
\subsection{$y2>y1$}
In this case, important is relation between two angular coefficients $k1$
and $k2$:\\
$k1=(A-y1)/(A-x1);\quad k2=(y2-y1)/(x2-x1)$.\\
$\bf \,\,\,k2\,>\,k1$.
 In this case the line $(p1,p2)$ crosses the upper side of square at
the point $(x31m,A)$ with $x31m=(A-y2)/k2+x2$, see Fig.1 , panel 1A.
\begin{figure} \includegraphics[scale=.65]{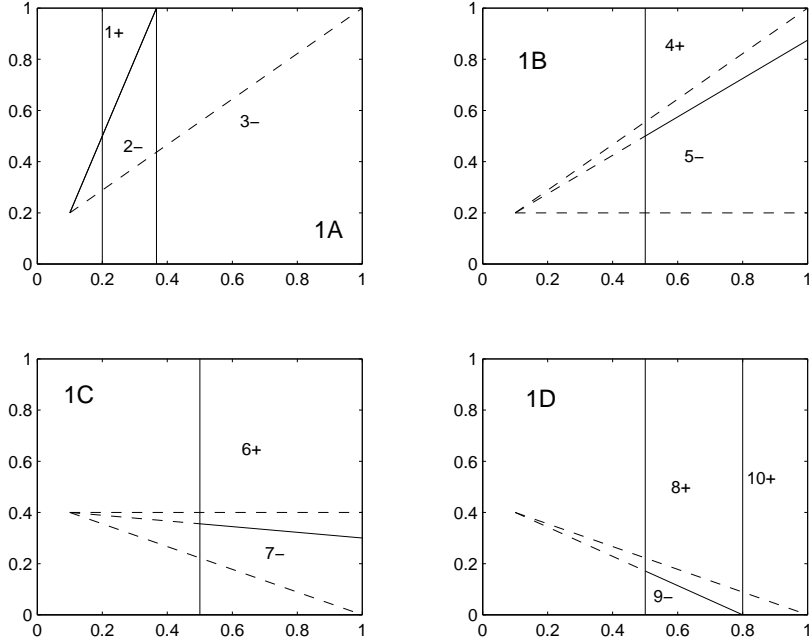}
\caption{Ten different cases of integration regions, bounded by solid
lines}\end{figure}
The region of integration over $x2,y2$ is:\\
$x1<x2<A;y21m<y2<A,y21m=k1(x2-x1)+y1.$\\
The first region of integration over $x3,y3$ is :
$x2<x3<x31m;\,\,k2(x3-x2)+y2<y3<A$, region 1+, panel 1A, Fig. 1.\\
At this region formula (1) gives positive value as points $(p1,p2,p3)$ make
right-handed system:
moving in direction $(p1 \rightarrow p2 \rightarrow p3 \rightarrow p1)$
we make counter-clockwise rotation.\\
Now we are ready to calculate first integral:
$$I1=\int_{0^{}}^A dx1\int_{0^{}}^A 
dy1\int_{x1^{}}^Adx2\int_{y21m^{}}^Ady2%
\int_{x2^{}}^{x31m}dx3\int_{k2(x3-x2)+y2^{}}^A(s)dy3=\frac{A^8}{34560}.$$
Second integral appears from region ''under'' the first integral's region, 
 region 2-, panel 1A, Fig. 1, and here formula (1) should be taken with
negative sign:
$$I2=\int_{0^{}}^Adx1\int_{0^{}}^Ady1\int_{x1^{}}^Adx2\int_{y21m^{}}^Ady2%
\int_{x2^{}}^{x31m}dx3\int_{0^{}}^{k2(x3-x2)+y2}(-s)dy3=\frac{23A^8}{34560}.$$
Note that 2nd integral differs from 1st one only by integration boundaries
over $y3$ (and by sign of $s$).\\
Also, the interesting exact relation occurs between numerical values of two
considered integrals : $I2=23I1.$\\
Last integral at the case $\bf k2>k1$, corresponding to region 3-, panel
1A, Fig. 1, is;  
$$I3=\int_{0^{}}^Adx1\int_{0^{}}^Ady1\int_{x1^{}}^Adx2\int_{y21m^{}}^Ady2%
\int_{x31m^{}}^Adx3\int_{0^{}}^A(-s)dy3=\frac{7A^8}{1728}=140I1.$$
Note that all three multiple integrals have the same first four particular
integral regions and we may write down them in a more compact form, but we
will not do this pure {''decorative''} operation.\\
 $\bf \,\,\,k2\,<\,k1$.
In this case and still at $\bf y1<y2$,  the line $(p1,p2)$ crosses the
right side of square;
now the region of integration over $y2$ is $y1<y2<y21m,$
and we have two different regions of integration over $p3$:\\
$x2<x3<A,$ $k2(x3-x2)+y2<y3<A,$ with positive $s$, region 4+, panel 1B,
Fig. 1,which gives $I4$, and\\
$x2<x3<A,0<y3<$ $k2(x3-x2)+y2,$ with negative $s$, r. 5-. p. 1A, Fig. 1,
which gives $I5$.\\
Therefore we have two additional integrals:
$$I4=\int_{0^{}}^Adx1\int_{0^{}}^Ady1\int_{x1^{}}^Adx2\int_{y1^{}}^{y21m}dy2%
\int_{x2^{}}^Adx3\int_{k2(x3-x2)+y2^{}}^A(s)dy3=\frac{19A^8}{34560}=19\,I1.$$
$$I5=\int_{0^{}}^Adx1\int_{0^{}}^Ady1\int_{x1^{}}^Adx2\int_{y1^{}}^{y21m}dy2%
\int_{x2^{}}^Adx3\int_{0^{}}^{k2(x3-x2)+y2^{}}(-s)dy3=\frac{37A^8}{34560}%
=37 I1.$$
\subsection{ $y2<y1$}
Now important is relation between two coefficients $k3$ and $k4:$\\
$k3=y1/(A-x1);\,\,k4=(y1-y2)/(x2-x1)$.\\
$\bf \,\,\,k4\,<\,k3$.
If $k4<k3,$ then the line $(p1,p2)$ crosses the right side of square ,
see panel 1C, Fig. 1, and we have the case completely
analogous to the case considered in the previous section (see also panel
1B) and two integrals, I6 with positive $s$ and I7 with
negative $s$, are equal to I4 and
I5 respectively. Here our geometrical approach is particularly explicitly
demonstrate his power:
it is sufficient to look at panels  1A - 1D of the Fig. 1 to be
convinced that actually we have only two different cases, one case 
when line (p1,p2) crosses two opposite sides of square and second case
when line (p1,p2) crosses two adjacent sides of square.\\
$\bf \,\,\,k4\,>\,k3$.
 Now the line (p1,p2) crosses lower side of square, panel 1D, Fig. 1,
and we have in essence the case coinciding with previous one, panel 1A,
Fig. 1, so we have another three integrals with values
actually found before: I8=I1, I9=I2, and I10=I3.\\
Now sum of all 10 integrals is equal to $II = 11A^8/864.$
The normalizing coefficient is found by calculating a sum of abovementioned
10 integrals with (+/-s)=1 in all
cases which gives $JJ=A^6/6,$ that is 1/6 of volume of hypercube with side
A. So the mean area of random triangle inside the square is
$II/JJ=11/144$ of host-figure's square, $A^2$.
\section{RT in rectangle}
Our geometrical approach allows easily to calculate also the mean area of
random triangle when host-figure is rectangle. Being experienced with 
the square case we consider here only two
cases leading to 5 integrals. \\
As result, we present a simple and transparent Mathematica's code for
calculating the mean area of random triangle in rectangle with sides A and
B.
\begin{verbatim}
(* y2\,>\,y1 *)
k1=(B-y1)/(A-x1);k2=(y2-y1)/(x2-x1); Y2=k1*(x2-x1)+y1; Y3=k2*(x3-x1)+y1;
(* k2\,>\,k1 *)
X=(B-y2)/k2+x2;
I1:=Integrate[s,x1,0,A,y1,0,B,x2,x1,A,y2,Y2,B, x3,x2,X,y3,Y3,B]; 
(* I1=A^4*B^4/34560 *)
I2:=Integrate[-s,x1,0,A,y1,0,B,x2,x1,A,y2,Y2,B, x3,x2,X,y3,0,Y3]; 
(* I2=23*I1 *)
I3:=Integrate[-s,x1,0,A,y1,0,B,x2,x1,A,y2,Y2,B, x3,X,A,y3,0,B]; 
(* I3=140 I1 *)
(* k2\,<\,k1 *)
I4:=Integrate[s,x1,0,A,y1,0,B,x2,x1,A,y2,y1,Y2, x3,x2,A,y3,Y3,B]; 
(* I4=19*I1 *)
I5:=Integrate[-s,x1,0,A,y1,0,B,x2,x1,A,y2,y1,Y2, x3,x2,A,y3,0,Y3]; 
(* I5=37*I1 *)
I15=I1+I2+I3+I4+I5; (* I15=11*A^4*B^4/1728 *)
(*calculation of normalizing coefficient*)
s=1;J1=I1;(* J1=A^3*B^3/432 *);s=-1;J2=I2;(* J2=5*J1 *)
s=-1;J3=I3;(* J3=18*J1 *);s=1;J4=I4;(* J4=5*J1 *)
s=-1;J5=I5;(* J5=7*J1 *); J15=J1+J2+J3+J4+J5; (*J15=A^3*B^3/12*);
MeanSquareOfRandomTriangleInRectangle=I15/J15;(* = (11/144)*A*B *)
(*MeanSquareOfRandomTriangleInRectangle = 11/144 of rectangle's square*)
 \end{verbatim}
\section{RT in square frame}
Now we consider the related problem of random triangle in square frame.
Let three points are randomly (with constant differential probability
function) distributed along the sides of unit square (side length and
square being 1). What is the mean area of triangles formed by these
points as vertices?\\
The solution is elementary, but we consider it for pure pedagogical purposes.
First observation is that due to symmetry of square and due to symmetry of
points it is sufficient
to assume that 1st particle is at bottom side of square.
Then four different cases should be considered.\\
{\bf 2nd particle is also at bottom side}\\
Let 3rd particle moves with constant
linear velocity along all sides
of the the square. We are looking for value of integral over coordinates of
3rd particle and then over
coordinates of 2nd particle which is allowed to move only along the bottom
side of square.\\
In the Mathematica's language we should calculate the following path
integral:
\begin{verbatim}
s[x1_,y1_,x2_,y2_,x3_,y3_]:=(1/2)*Abs[(x1*(y2-y3)+
x2*(-y1+y3)+x3*(y1-y2))]:
I1:=Integrate[(Integrate[s[x1,0,x2,0,x3,0],{ x3,0,1}]+
Integrate[s[x1,0,x2,0,1,y3],{ x2, 0,1},{y 3,0,1}]+
Integrate[s[x1,0,x2,0,x3,1],{ x2, 0,1},{ x3,0,1}]+
Integrate[s[x1,0,x2,0,0,y3],{ x2, 0,1},{ y3,0,1}]),{x2,0,1}]
\end{verbatim}
Result is $I1=\frac 12-x1+x1^2.$\\
{\bf 2nd particle is at the right side of the square}\\
The relevant integral which we do not write down
is equal to $I2=\frac{11-8x1+3x1^2}{12}.$\\
{\bf 2nd particle is at the upper side of the square}.
$I3=\frac{11-6x1+6x1^2}{12}.$\\
{\bf 2nd particle is at the left side of the square}.
$I4=\frac{6+2x1+3x1^2}{12}.$\\
Sum of these for integrals gives $I14=I1+I2+I3+I4=\frac{17}6-2x1+2x1^2.$
Now \verb(Integrate[I14,{x1,0,1}]( 
gives 5/2. The normalizing coefficient is evidently 16, as we calculate 4
path integrals (over 3rd particle) which of them has path length equal to 4.
Final result is: the mean area of random triangle inscribed in unit square
is 5/32.\\
Let us look at this value (and check it!) from another point of view. We
divide all sides of unit square to 10 equal parts and
let each of three particles take all mid-points of these 40 parts. Then we
have 40x40x40x40=640,000 triangles with mean area (as calculated by
Mathematica) equal to 249/1600=498/3200 that is very close to 5/32.\\
So we understand more vividly in what sense the mean area of random
triangles inscribed in the square is 5/32.\\
It is very interesting to compare these two ''mean" values 11/144 and
5/32. First value is 22/45, that is almost exactly 1/2, of the second one.
That is mean square of triangles with vertices randomly distributed
all over the square is almost exactly half of mean square of triangles
vertices of which are allowed to occur only at sides of square. \\
The reason of considering this last problem is originally related to my
attempts to find simple solution of first problem. Is seemed to me
that by solving the problem of inscribed triangles I could somehow  solve
the first problem also. Some hazy ideas about differentiation/integration
connection between two problems unfortunately gave no yield and I solved
these two problems separately.\\
What is left is the problem of mean volume of tetrahedron in the cube 
(M.Trott, personal communication ). I hope that geometrical approach will
also help in this  much more difficult problem. But if the geometrical
approach managed to reduce the number of integrals from 496 cumbersome
ones in original solution by M. Trott to 5 very simple integrals, 
hopefully it will help in tetrahedron-in-cube problem as well. \\
Numerical value get by Mathematica gives 1/72, but this is not
exact value, this is value get by Mathematica's command 
{\verb(Rationalize[NumericalValue,10^-4](}.
\section*{Acknowledgement}
The useful correspondence with M. Trott is highly appreciated.
\section*{References}
[1] M. Trott, Mathematica Journal, v7 i2, 189-197,1998.
\end{document}